\documentclass{article}
\usepackage{amsfonts}
\usepackage{amsmath}

\input{tcilatex}

\begin{document}

\begin{center}
{\Huge Derivative with two fractional orders: New door of investigation
toward revolution in fractional calculus}$\bigskip $

Abdon Atangana

Institute for Groundwater studies, Faculty of Natural and Agricultural
Sciences, University of the Free State, 9300, South Africa,

email AtanganaA@ufs.ac.za
\end{center}

\bigskip

\noindent \noindent \textbf{Abstract: }In order to describe more complex
problem using the concept of fractional derivative, we introduce in this
paper a derivative with two fractional orders. Each order therefore plays an
important rule while modeling for instance problems with two layers with
different properties. This is the case for instance in thermal science for a
reaction diffusion within a media with two different layers with different
properties. Another case of this study can be used in the study of
groundwater flowing with an aquifer where the ground has layers with
different properties. The paper presents some basic definitions and
properties.

\noindent \noindent \textbf{Keywords}: New concept of fractional derivative,
bi-order, generalized Mittag-Leffler function.

\section{Introduction}

The concept of fractional calculus was introduced by mathematicians to solve
problems that could not be handled by local derivative, the fractional order
alpha that appears in the concept of fractional derivative can be used to
represent some physical parameters. Nevertheless it is not possible for the
fractional derivative to be used in order to describe the movement of for
instance heat via material with different layers where each layer possesses
different material from each other. Riemann-Liouville introduced the
fractional derivative as a derivative with convolution of a given function
and the power law function [1-5]. Caputo introduced the fractional derivative as
convolution of derivative of a given function with power law function.
Caputo and Fabrizio introduced the fractional derivative as convolution of
derivative of a given function with exponential [2-4]. Atangana and Goufo
independently with Caputo and Fabrizio introduced another one as derivative
of a convolution of a given function with the exponential function. Atangana
and Dumitru introduced fractional derivative as a both convolution of a
derivative of a given function and the Mittag-Leffler function [4-5]. These
concepts of fractional derivative so far has been used with great success
and misused also. Nevertheless they are not able to handle the concept of
heterogeneity with great success. The concept of Brownian motion was
introduced using fractional derivative to describe to choice that a particle
has when moving via heterogeneous media, however, this concept failed due to
the fact that, the particle does not chose where to move within a given
heterogeneous media, but it condition by the properties of the media [6-7]. For
instance, the drop of water is moving via a given geological formation under
the influence of the properties of the soul. It is therefore important to
introduce a suitable derivative that can used to handle this type of
physical problem. In this paper, we introduce the concept of fractional
order with two orders alpha and beta.

\section{New concept of fractional derivatives}

In this section, some new definitions of fractional derivatives with
bi-order are presented.

\noindent \noindent \textbf{Definition 1 }: Let $f(x)$ defined in $[a,X]$
necessary differentiable such that for all $x\in \lbrack a,X]$ the
convolution of $x^{1-\alpha -n}E_{\beta }\left\{ -\frac{\beta }{1-\beta }%
x^{\beta +\alpha }\right\} $ with the function $\frac{d^{n}f(x)}{dx^{n}}$ exists, then the
fractional derivative of the function $f$ with order $\alpha ,\beta \in
(0,1) $ is given as:%
\begin{eqnarray}
&&_{a}^{AC}D_{x}^{\alpha ,\beta }f(x) \\
&=&\frac{A(\beta )}{1-\beta }\frac{1}{\Gamma \left\{ 1-\alpha \right\} }%
\int\limits_{a}^{x}\frac{d^{n}f(t)}{dt^{n}}(x-t)^{1-\alpha -n}E_{\beta
}\left\{ -\frac{\beta }{1-\beta }(x-t)^{\beta +\alpha }\right\} dt,  \notag
\\
0 &<&\beta <1,\text{ \ \ }n-1<\alpha <n.  \notag
\end{eqnarray}

\noindent \noindent \textbf{Definition 2 }: Let $f(x)$ defined in $[a,X]$
not necessary differentiable such that for all $x\in \lbrack a,X]$ the
convolution of $x^{1-\alpha -n}E_{\beta }\left\{ -\frac{\beta }{1-\beta }%
x^{\beta +\alpha }\right\} $ with the function $f(x)$ exists, then the
fractional derivative of the function $f$ with order $\alpha ,\beta \in
(0,1) $ is given as:%
\begin{eqnarray}
&&_{a}^{AR}D_{x}^{\alpha ,\beta }f(x) \\
&=&\frac{A(\beta )}{1-\beta }\frac{1}{\Gamma \left\{ 1-\alpha \right\} }%
\frac{d^{n}}{dx^{n}}\int\limits_{a}^{x}f(t)(x-t)^{1-\alpha -n}E_{\beta
}\left\{ -\frac{\beta }{1-\beta }(x-t)^{\beta +\alpha }\right\} dt,  \notag
\\
0 &<&\beta <1,\text{ \ \ }n-1<\alpha <n.  \notag
\end{eqnarray}%
\noindent \noindent \textbf{Definition 3 }: Let $f(x)$ defined in $[X,\infty
]$ not necessary differentiable such that for all $x\in \lbrack 0,\infty ]$
the convolution of $x^{1-\alpha -n}E_{\beta }\left\{ -\frac{\beta }{1-\beta }%
x^{\beta +\alpha }\right\} $ with the function $\frac{d^{n}f(x)}{dx^{n}}$ exists, then the
fractional derivative of the function $f$ with order $\alpha ,\beta \in
(0,1) $ is given as:%
\begin{eqnarray}
&&_{a}^{AR}D_{x}^{\alpha ,\beta }f(x)  \notag \\
&=&\frac{A(\beta )}{1-\beta }\frac{1}{\Gamma \left\{ 1-\alpha \right\} }%
\frac{d^{n}}{dx^{n}}\int\limits_{x}^{\infty }f(t)(t-x)^{1-\alpha -n}E_{\beta
}\left\{ -\frac{\beta }{1-\beta }(x-t)^{\beta +\alpha }\right\} dt,  \notag
\\
0 &<&\beta <1,\text{ \ \ }n-1<\alpha <n.
\end{eqnarray}%
\noindent \noindent \textbf{Definition 4 }: Let $f(x)$ defined in $[X,\infty
]$ necessary differentiable such that for all $x\in \lbrack 0,\infty ]$ the
convolution of $x^{1-\alpha -n}E_{\beta }\left\{ -\frac{\beta }{1-\beta }%
x^{\beta +\alpha }\right\} $ with the function $f(x)$ exists, then the
fractional derivative of the function $f$ with order $\alpha ,\beta \in
(0,1) $ is given as:%
\begin{eqnarray}
&&_{a}^{AC}D_{x}^{\alpha ,\beta }f(x) \\
&=&\frac{A(\beta )}{1-\beta }\frac{1}{\Gamma \left\{ 1-\alpha \right\} }%
\int\limits_{x}^{\infty }\frac{d^{n}f(t)}{dt^{n}}(t-x)^{1-\alpha -n}E_{\beta
}\left\{ -\frac{\beta }{1-\beta }(x-t)^{\beta +\alpha }\right\} dt,  \notag
\\
0 &<&\beta <1,\text{ \ \ }n-1<\alpha <n.  \notag
\end{eqnarray}%
\noindent \noindent \textbf{Definition 5 }: Let $f(x)$ defined in $[-\infty
,X]$ not necessary differentiable such that for all $x\in \lbrack -\infty
,X] $ the convolution of $x^{1-\alpha -n}E_{\beta }\left\{ -\frac{\beta }{%
1-\beta }x^{\beta +\alpha }\right\} $ with the function $f(x)$ exists, then
the fractional derivative of the function $f$ with order $\alpha ,\beta \in
(0,1)$ is given as:%
\begin{eqnarray*}
&&_{a}^{AR}D_{x}^{\alpha ,\beta }f(x) \\
&=&\frac{A(\beta )}{1-\beta }\frac{1}{\Gamma \left\{ 1-\alpha \right\} }%
\frac{d^{n}}{dx^{n}}\int\limits_{-\infty }^{x}f(t)(x-t)^{1-\alpha
-n}E_{\beta }\left\{ -\frac{\beta }{1-\beta }(x-t)^{\beta +\alpha }\right\}
dt, \\
0 &<&\beta <1,\text{ \ \ }n-1<\alpha <n.
\end{eqnarray*}%
\noindent \noindent \textbf{Definition 6 }: Let $f(x)$ defined in $[-\infty
,X]$ necessary differentiable such that for all $x\in \lbrack -\infty ,X]$
the convolution of $x^{1-\alpha -n}E_{\beta }\left\{ -\frac{\beta }{1-\beta }%
x^{\beta +\alpha }\right\} $ with the function $f(x)$ exists, then the
fractional derivative of the function $f$ with order $\alpha ,\beta \in
(0,1) $ is given as:%
\begin{eqnarray}
&&_{a}^{AC}D_{x}^{\alpha ,\beta }f(x) \\
&=&\frac{A(\beta )}{1-\beta }\frac{1}{\Gamma \left\{ 1-\alpha \right\} }%
\int\limits_{-\infty }^{x}\frac{d^{n}f(t)}{dt^{n}}(x-t)^{1-\alpha
-n}E_{\beta }\left\{ -\frac{\beta }{1-\beta }(x-t)^{\beta +\alpha }\right\}
dt,  \notag \\
0 &<&\beta <1,\text{ \ \ }n-1<\alpha <n.  \notag
\end{eqnarray}

\section{New properties of the A.R and A.C derivatives}

Some properties of the A.R and A.C derivatives with two order $\alpha ,\beta 
$ are given below:

\noindent \noindent \textbf{The Laplace transform of A.R:}

The Laplace transform of A.R derivative is obtained via Laplace transform
definition below: 
\begin{eqnarray}
&&L\left\{ _{a}^{AR}D_{x}^{\alpha ,\beta }f(x)\right\} \\
&=&\frac{A(\beta )}{1-\beta }\frac{1}{\Gamma \left\{ 1-\alpha \right\} }%
\int\limits_{0}^{\infty }e^{-pt}\frac{d}{dt}\int\limits_{0}^{t}f(\tau
)(t-\tau )^{-\alpha }E_{\beta }\left\{ -\frac{\beta }{1-\beta }(t-\tau
)^{\beta +\alpha }\right\} d\tau dt  \notag
\end{eqnarray}%
\begin{eqnarray}
&=&\frac{A(\beta )}{1-\beta }\frac{p}{\Gamma \left\{ 1-\alpha \right\} }%
L\left\{ t^{-\alpha }E_{\beta }\left( -\frac{\beta }{1-\beta }t^{\beta
+\alpha }\right) \right\} L\left\{ f(t)\right\} \\
&=&\frac{A(\beta )}{1-\beta }\frac{p}{\Gamma \left\{ 1-\alpha \right\} }%
L\left\{ t^{1-\alpha -1}E_{\beta }\left( -\frac{\beta }{1-\beta }t^{\beta
+\alpha }\right) \right\} F(p)  \notag \\
&=&\frac{A(\beta )}{1-\beta }\frac{p}{\Gamma \left\{ 1-\alpha \right\} }%
L\left\{ t^{p-1}E_{\beta }\left( -\frac{\beta }{1-\beta }t^{\sigma }\right)
\right\} F(p)  \notag
\end{eqnarray}%
However we have from the literature that 
\begin{equation}
L\left\{ t^{p-1}E_{\beta ,\Psi }^{\delta }\left( at^{\alpha }\right)
\right\} =\frac{s^{-p}}{\Gamma (\delta )}\text{ }_{2}\Psi _{1}\left[ 
\begin{array}{c}
(\delta ,1),(\Psi ,\alpha ); \\ 
(\gamma ,\beta )\text{ \ \ \ \ \ \ \ ;}%
\end{array}%
\begin{array}{c}
\frac{a}{s^{\alpha }}%
\end{array}%
\right]
\end{equation}%
for $\alpha ,p,\delta ,\beta ,\gamma >0$ thus if we take $\delta =\Psi =1,$ $%
a=-\frac{\beta }{1-\beta },$ and $\sigma =\alpha +\beta $ then we obtain the
following

\begin{eqnarray}
&&L\left\{ t^{p-1}E_{\beta }\left( -\frac{\beta }{1-\beta }t^{\sigma
}\right) \right\}  \notag \\
&=&s_{\text{ \ \ \ \ \ \ }2}^{\alpha -1}\Psi _{1}\left[ 
\begin{array}{c}
(1,1),(1-\alpha ,\alpha +\beta ); \\ 
(1,\beta )\text{ \ \ \ \ \ \ \ \ \ \ \ \ \ \ \ \ \ \ ;}%
\end{array}%
\begin{array}{c}
-\frac{\beta }{1-\beta }.\frac{1}{s^{\alpha +\beta }}%
\end{array}%
\right]
\end{eqnarray}%
Where the function $_{p}\Psi _{q}$ is the wright's generalized
hyper-geometric function defined by means of the similar representation in
the following form [8-9]:%
\begin{equation}
_{p}\Psi _{q}(z)=\sum\limits_{r=0}^{\infty }\frac{\left\{
\prod\limits_{j=1}^{p}r(a_{j}+A_{j}r)\right\} }{\left\{
\prod\limits_{j=1}^{q}r(b_{j}+B_{j}r)\right\} }\frac{z^{r}}{r!},
\end{equation}%
therefore replacing the above we obtain%
\begin{eqnarray}
&&L\left\{ _{a}^{AR}D_{t}^{\alpha ,\beta }f(t)\right\} \\
&=&\frac{A(\beta )}{1-\beta }\frac{1}{\Gamma \left\{ 1-\alpha \right\} }s_{%
\text{ \ \ \ \ \ \ }2}^{\alpha }\Psi _{1}\left[ 
\begin{array}{c}
(1,1),(1,\alpha +\beta ); \\ 
(1,\beta )\text{ \ \ \ \ \ \ \ \ \ \ \ \ \ ;}%
\end{array}%
\begin{array}{c}
-\frac{\beta }{1-\beta }.\frac{1}{s^{\alpha +\beta }}%
\end{array}%
\right]  \notag
\end{eqnarray}

\noindent \noindent \textbf{Existence of A.C:}

\noindent \noindent \textbf{Theorem 1 }: Let $f$ be continuous in $[0,T]$
and $\alpha ,\beta \in (0,1)$ then the following inequality is obtained if 
\begin{equation}
\left\Vert \frac{df(t)}{dt}\right\Vert \leq \Theta \left\Vert
f(t)\right\Vert _{C}
\end{equation}

\noindent \noindent \textbf{Proof} :%
\begin{eqnarray}
&&\left\Vert _{0}^{AC}D_{t}^{\alpha ,\beta }f(t)\right\Vert _{C} \\
&=&\frac{A(\beta )}{1-\beta }\frac{1}{\Gamma \left\{ 1-\alpha \right\} }%
\left\Vert \int\limits_{0}^{t}\frac{df(\tau )}{d\tau }(t-\tau )^{-\alpha
}E_{\beta }\left\{ -\frac{\beta }{1-\beta }(t-\tau )^{\beta +\alpha
}\right\} d\tau \right\Vert  \notag \\
&\leq &\frac{A(\beta )}{1-\beta }\frac{1}{\Gamma \left\{ 1-\alpha \right\} }%
\int\limits_{0}^{t}\left\Vert \frac{df(\tau )}{d\tau }\right\Vert \left\Vert
(t-\tau )^{-\alpha }E_{\beta }\left\{ -\frac{\beta }{1-\beta }(t-\tau
)^{\beta +\alpha }\right\} \right\Vert d\tau  \notag \\
&\leq &\frac{A(\beta )}{1-\beta }\frac{1}{\Gamma \left\{ 1-\alpha \right\} }%
\Theta \left\Vert f(t)\right\Vert _{C}\int\limits_{0}^{t}(t-\tau )^{-\alpha
}E_{\beta }\left\{ -\frac{\beta }{1-\beta }(t-\tau )^{\beta +\alpha
}\right\} d\tau  \notag \\
&\leq &\frac{A(\beta )}{1-\beta }\frac{1}{\Gamma \left\{ 1-\alpha \right\} }%
\Theta \left\Vert f\right\Vert _{C}\int\limits_{0}^{t}(t-\tau )^{1-\alpha
-1}E_{\beta }\left\{ -\frac{\beta }{1-\beta }(t-\tau )^{\beta +\alpha
}\right\} d\tau  \notag \\
&\leq &\frac{A(\beta )}{1-\beta }\frac{1}{\Gamma \left\{ 1-\alpha \right\} }%
\Theta \left\Vert f\right\Vert _{C}\int\limits_{0}^{t}(t-\tau
)^{p-1}E_{\beta }\left\{ -\frac{\beta }{1-\beta }(t-\tau )^{\beta +\alpha
}\right\} d\tau .  \notag
\end{eqnarray}%
Let $y=t-\tau $ and $p=1-\alpha $ then 
\begin{equation}
\int\limits_{0}^{t}(t-\tau )^{p-1}E_{\beta }\left\{ -\frac{\beta }{1-\beta }%
(t-\tau )^{\beta +\alpha }\right\} d\tau =\int\limits_{0}^{t}y^{p-1}E_{\beta
}\left\{ -\frac{\beta }{1-\beta }y^{\beta +\alpha }\right\} dy.
\end{equation}%
However%
\begin{equation}
\int\limits_{0}^{t}y^{p-1}E_{\beta }\left\{ -\frac{\beta }{1-\beta }y^{\beta
+\alpha }\right\} dy=t^{p}E_{\beta ,2}\left\{ -\frac{\beta }{1-\beta }%
t^{\beta +\alpha }\right\}
\end{equation}%
then 
\begin{eqnarray}
\left\Vert _{0}^{AC}D_{t}^{\alpha ,\beta }f(t)\right\Vert _{C} &\leq
&\left\Vert f\right\Vert _{C}K\text{ ,} \\
K &=&\frac{A(\beta )}{1-\beta }\frac{1}{\Gamma \left\{ 1-\alpha \right\} }%
\Theta T^{1-\alpha }E_{\beta ,2}\left\{ -\frac{\beta }{1-\beta }T^{\beta
+\alpha }\right\}  \notag
\end{eqnarray}

\noindent \noindent \textbf{Test of the Lipschitz continuity:}

Let $f$ and $g$ be two defined function on $[0,T]$ such that their A.C
derivatives exist then%
\begin{eqnarray}
&&\left\Vert _{0}^{AC}D_{t}^{\alpha ,\beta }f(t)-_{0}^{AC}D_{t}^{\alpha
,\beta }g(t)\right\Vert _{C} \\
&=&\left\Vert \frac{A(\beta )}{1-\beta }\frac{1}{\Gamma \left\{ 1-\alpha
\right\} }\int\limits_{0}^{t}\frac{d(f-g)(\tau )}{d\tau }(t-\tau )^{-\alpha
}E_{\beta }\left\{ -\frac{\beta }{1-\beta }(t-\tau )^{\beta +\alpha
}\right\} d\tau \right\Vert  \notag \\
&\leq &\frac{A(\beta )}{1-\beta }\frac{1}{\Gamma \left\{ 1-\alpha \right\} }%
\int\limits_{0}^{t}\Theta _{1}\left\Vert f-g\right\Vert _{C}\left\Vert
(t-\tau )^{-\alpha }E_{\beta }\left\{ -\frac{\beta }{1-\beta }(t-\tau
)^{\beta +\alpha }\right\} \right\Vert d\tau  \notag \\
&\leq &\frac{A(\beta )}{1-\beta }\frac{\Theta _{1}}{\Gamma \left\{ 1-\alpha
\right\} }\left\Vert f-g\right\Vert _{C}\int\limits_{0}^{t}(t-\tau
)^{-\alpha }E_{\beta }\left\{ -\frac{\beta }{1-\beta }(t-\tau )^{\beta
+\alpha }\right\} d\tau  \notag \\
&\leq &\frac{A(\beta )}{1-\beta }\frac{\Theta _{1}}{\Gamma \left\{ 1-\alpha
\right\} }T^{-\alpha }E_{\beta }\left\{ -\frac{\beta }{1-\beta }T^{\beta
+\alpha }\right\} \left\Vert f-g\right\Vert _{C}  \notag \\
&\leq &L\left\Vert f-g\right\Vert _{C}  \notag
\end{eqnarray}%
Thus 
\begin{equation}
\left\Vert _{0}^{AC}D_{t}^{\alpha ,\beta }f(t)-_{0}^{AC}D_{t}^{\alpha ,\beta
}g(t)\right\Vert _{C}\leq L\left\Vert f(t)-g(t)\right\Vert _{C}.
\end{equation}%
This implies A.C has Lipschitz continuity condition.

\noindent \noindent \textbf{Theorem 2 }: Let $f$ be a continuous function
defined in $[0,T]$ such that $_{0}^{AC}D_{t}^{\alpha ,\beta }f(t)$ and $%
_{0}^{AR}D_{t}^{\alpha ,\beta }f(t)$ exists, then the following relation is
established.%
\begin{equation}
_{0}^{AC}D_{t}^{\alpha ,\beta }f(t)=-G(t)+_{0}^{AR}D_{t}^{\alpha ,\beta
}f(t).
\end{equation}

\noindent \noindent Proof : We achieve this using the Laplace transform%
\begin{eqnarray}
&&L\left\{ _{0}^{AC}D_{t}^{\alpha ,\beta }f(t)\right\} \\
&=&L\left\{ \frac{A(\beta )}{1-\beta }\frac{1}{\Gamma \left\{ 1-\alpha
\right\} }\int\limits_{0}^{t}\frac{df(\tau )}{d\tau }(t-\tau )^{-\alpha
}E_{\beta }\left\{ -\frac{\beta }{1-\beta }(t-\tau )^{\beta +\alpha
}\right\} d\tau \right\}  \notag \\
&=&\frac{A(\beta )}{1-\beta }\frac{1}{\Gamma \left\{ 1-\alpha \right\} }%
L\left\{ \frac{df(\tau )}{d\tau }\right\} L\left\{ t^{-\alpha }E_{\beta
}\left\{ -\frac{\beta }{1-\beta }t^{\beta +\alpha }\right\} \right\}  \notag
\\
&=&\frac{A(\beta )}{1-\beta }\frac{1}{\Gamma \left\{ 1-\alpha \right\} }%
\left\{ pF(p)-f(0)\right\} L\left\{ t^{-\alpha }E_{\beta }\left\{ -\frac{%
\beta }{1-\beta }t^{\beta +\alpha }\right\} \right\}  \notag \\
&=&\frac{A(\beta )}{1-\beta }\frac{p}{\Gamma \left\{ 1-\alpha \right\} }%
\text{ }_{2}\Psi _{1}\left[ 
\begin{array}{c}
(1,1),(1-\alpha ,\alpha +\beta ); \\ 
(1,\beta )\text{ \ \ \ \ \ \ \ \ \ \ \ \ \ \ \ \ \ \ ;}%
\end{array}%
\begin{array}{c}
-\frac{\beta }{1-\beta }.\frac{1}{s^{\alpha +\beta }}%
\end{array}%
\right] F(p)  \notag \\
&&-\frac{A(\beta )}{1-\beta }\frac{p^{\alpha -1}}{\Gamma \left\{ 1-\alpha
\right\} }_{2}\Psi _{1}\left[ 
\begin{array}{c}
(1,1),(1-\alpha ,\alpha +\beta ); \\ 
(1,\beta )\text{ \ \ \ \ \ \ \ \ \ \ \ \ \ \ \ \ \ \ ;}%
\end{array}%
\begin{array}{c}
-\frac{\beta }{1-\beta }.\frac{1}{s^{\alpha +\beta }}%
\end{array}%
\right] f(0).  \notag
\end{eqnarray}

Therefore the inverse Laplace transform applied on the above yields,%
\begin{eqnarray}
_{0}^{AC}D_{t}^{\alpha ,\beta }f(t) &=&_{0}^{AR}D_{t}^{\alpha ,\beta }f(t)-%
\frac{A(\beta )}{1-\beta }\frac{1}{\Gamma \left\{ 1-\alpha \right\} }%
f(0)t^{-\alpha }E_{\beta }\left\{ -t^{\beta +\alpha }\right\} \\
&=&_{0}^{AR}D_{t}^{\alpha ,\beta }f(t)-G(t).  \notag
\end{eqnarray}%
This completes the proof.

\noindent \noindent \textbf{Sumudu Transform of A.C:}

The Sumudu transform of A.C is expressed as follows:%
\begin{eqnarray}
&&S\left( _{0}^{AR}D_{t}^{\alpha ,\beta }f(t)\right) \\
&=&\frac{A(\beta )}{1-\beta }\frac{1}{\Gamma \left\{ 1-\alpha \right\} }%
S\left( \frac{d}{dt}\int\limits_{0}^{t}f(\tau )(t-\tau )^{-\alpha }E_{\beta
}\left\{ -\frac{\beta }{1-\beta }(t-\tau )^{\beta +\alpha }\right\} d\tau
\right)  \notag \\
&=&\frac{A(\beta )}{1-\beta }\frac{1}{\Gamma \left\{ 1-\alpha \right\} }%
\frac{1}{u}.S\left( \int\limits_{0}^{t}f(\tau )(t-\tau )^{-\alpha }E_{\beta
}\left\{ -\frac{\beta }{1-\beta }(t-\tau )^{\beta +\alpha }\right\} d\tau
\right)  \notag \\
&=&\frac{A(\beta )}{1-\beta }\frac{1}{\Gamma \left\{ 1-\alpha \right\} }%
\frac{1}{u}S\left( f(t)\right) .S\left( t^{-\alpha }E_{\beta }\left\{ -\frac{%
\beta }{1-\beta }t^{\beta +\alpha }\right\} \right)  \notag
\end{eqnarray}%
where 
\begin{eqnarray}
S\left( t^{-\alpha }E_{\beta }\left\{ -\frac{\beta }{1-\beta }t^{\beta
+\alpha }\right\} \right) &=&\int\limits_{0}^{\infty }e^{-t}u^{-\alpha
}t^{-\alpha }E_{\beta }\left\{ -\frac{\beta }{1-\beta }(ut)^{\beta +\alpha
}\right\} dt \\
&=&u^{-\alpha }\int\limits_{0}^{\infty }e^{-t}t^{-\alpha }E_{\beta }\left\{
-\lambda t^{\beta +\alpha }\right\} dt,  \notag
\end{eqnarray}%
where $\lambda =\frac{\beta }{1-\beta }u^{\alpha +\beta }.$%
\begin{eqnarray}
S\left( t^{-\alpha }E_{\beta }\left\{ -\frac{\beta }{1-\beta }t^{\beta
+\alpha }\right\} \right) &=&u^{-\alpha }\int\limits_{0}^{\infty
}e^{-t}t^{-\alpha }E_{\beta }\left\{ -\lambda t^{\beta +\alpha }\right\} dt
\\
&=&u^{-\alpha }\int\limits_{0}^{\infty }e^{-t}t^{-\alpha
}\sum\limits_{j=0}^{\infty }\frac{\left( -\lambda t^{\beta +\alpha }\right)
^{j}}{\Gamma \left( \alpha j+1\right) }  \notag
\end{eqnarray}%
where $\lambda =\frac{\beta }{1-\beta }u^{\alpha +\beta }.$%
\begin{eqnarray}
&=&u^{-\alpha }\sum\limits_{j=0}^{\infty }\frac{\left( -\lambda \right) ^{j}%
}{\Gamma \left( \alpha j+1\right) }\int\limits_{0}^{\infty }e^{-t}t^{\alpha
j+\beta j-\alpha }dt \\
&=&u^{-\alpha }\sum\limits_{j=0}^{\infty }\frac{\left( -\lambda \right) ^{j}%
}{\Gamma \left( \alpha j+1\right) }\int\limits_{0}^{\infty }e^{-t}t^{\alpha
j+\beta j-\alpha +1-1}dt  \notag \\
&=&u^{-\alpha }\sum\limits_{j=0}^{\infty }\frac{\left( -\lambda \right) ^{j}%
}{\Gamma \left( \alpha j+1\right) }\int\limits_{0}^{\infty }e^{-t}t^{\Theta
-1}dt  \notag
\end{eqnarray}%
where $\Theta =\alpha j+\beta j-\alpha +1.$%
\begin{eqnarray}
&=&u^{-\alpha }\sum\limits_{j=0}^{\infty }\frac{\left( -\lambda \right) ^{j}%
}{\Gamma \left( \alpha j+1\right) }\Gamma (\Theta ) \\
&=&u^{-\alpha }\sum\limits_{j=0}^{\infty }\frac{\left( -\lambda \right) ^{j}%
}{\Gamma \left( \alpha j+1\right) }\Gamma (\alpha j+\beta j-\alpha +1) 
\notag \\
&=&u^{-\alpha }\sum\limits_{j=0}^{\infty }\frac{\left( -\frac{\beta }{%
1-\beta }u^{\beta +\alpha }\right) ^{j}\Gamma (\alpha j+\beta j-\alpha +1)}{%
\Gamma \left( \alpha j+1\right) }  \notag
\end{eqnarray}%
Therefore the Sumudu transform of $_{0}^{AR}D_{t}^{\alpha ,\beta }f(t)$ is
given as 
\begin{eqnarray}
&&S\left( _{0}^{AR}D_{t}^{\alpha ,\beta }f(t)\right) \\
&=&\frac{A(\beta )}{1-\beta }\frac{1}{\Gamma \left\{ 1-\alpha \right\} }%
u^{-\alpha -1}\sum\limits_{j=0}^{\infty }\frac{\left( -\frac{\beta }{1-\beta 
}u^{\beta +\alpha }\right) ^{j}\Gamma (\alpha j+\beta j-\alpha +1)}{\Gamma
\left( \alpha j+1\right) }  \notag
\end{eqnarray}

\section{Partial derivative with A.C and A.R}

Since one can model many natural occurrence using the space and time, it is
therefore important \ to propose in this section the partial A.C and A.R
fractional derivative with order $\alpha ,\beta .$

\noindent \noindent \textbf{Definition 7 }: Let $f(x,t)$ be a function
differentiable in $x$ or $t$ direction. Let $\alpha ,\beta \in (0,1)$ such
that the convolution of $x^{-\alpha }E_{\beta }\left\{ -\frac{\beta }{%
1-\beta }x^{\beta +\alpha }\right\} $ and $\frac{\partial f}{\partial x}$
exists than the A.C partial fractional derivative of $f$ of order $\alpha
,\beta $ is given as:%
\begin{eqnarray}
&&_{0}^{AC}D_{x}^{\alpha ,\beta }f(x,t) \\
&=&\frac{A(\beta )}{1-\beta }\frac{1}{\Gamma \left\{ 1-\alpha \right\} }%
\int\limits_{0}^{x}\frac{\partial f(\xi ,t)}{\partial t^{n}}(x-\xi
)^{-\alpha }E_{\beta }\left\{ -\frac{\beta }{1-\beta }(x-\xi )^{\beta
+\alpha }\right\} d\xi .  \notag
\end{eqnarray}%
\noindent \noindent \textbf{Definition 8 }: Let $f(x,t)$ be a function not
necessary differentiable in $x$ or $t$ direction. Let $\alpha ,\beta \in
(0,1)$ such that the convolution of $x^{-\alpha }E_{\beta }\left\{ -\frac{%
\beta }{1-\beta }x^{\beta +\alpha }\right\} $ and $\frac{\partial f}{%
\partial x}$ exists than the A.R partial fractional derivative of $f$ of
order $\alpha ,\beta $ is given as:%
\begin{eqnarray}
&&_{0}^{AR}D_{x}^{\alpha ,\beta }f(x,t) \\
&=&\frac{A(\beta )}{1-\beta }\frac{1}{\Gamma \left\{ 1-\alpha \right\} }%
\frac{\partial }{\partial x}\int\limits_{0}^{x}f(\xi ,t)(x-\xi )^{-\alpha
}E_{\beta }\left\{ -\frac{\beta }{1-\beta }(x-\xi )^{\beta +\alpha }\right\}
d\xi .  \notag
\end{eqnarray}%
\noindent \noindent \textbf{Definition 9 }: Let $f(x,t)$ be a function such
that $\frac{\partial ^{2}f}{\partial x\partial t}$ exists. Let $\alpha
,\beta \in (0,1)$ then A.C partial derivative of $f$ of order $\alpha ,\beta 
$ is given as:%
\begin{eqnarray}
&&_{0}^{AC}D_{x,t}^{\alpha ,\beta }f(x,t) \\
&=&\frac{A(\beta )}{1-\beta }\frac{1}{\Gamma \left\{ 1-\alpha \right\} }%
\int\limits_{0}^{x}\int\limits_{0}^{t}\frac{\partial ^{2}f(\xi ,\lambda )}{%
\partial \xi \partial \lambda }(x-\xi )^{-\alpha }(t-\lambda )^{-\alpha } 
\notag \\
&&.E_{\beta }\left\{ -\frac{\beta }{1-\beta }(x-\xi )^{\beta +\alpha
}\right\} E_{\beta }\left\{ -\frac{\beta }{1-\beta }(t-\lambda )^{\beta
+\alpha }\right\} d\xi d\lambda .  \notag
\end{eqnarray}%
However if $\frac{\partial ^{2}f}{\partial x\partial t}$ does not exist then 
\begin{eqnarray}
&&_{0}^{AR}D_{x,t}^{\alpha ,\beta }f(x,t) \\
&=&\frac{A(\beta )}{1-\beta }\frac{1}{\Gamma \left\{ 1-\alpha \right\} }%
\frac{\partial ^{2}f}{\partial x\partial t}\int\limits_{0}^{x}\int%
\limits_{0}^{t}f(\xi ,\lambda )(x-\xi )^{-\alpha }(t-\lambda )^{-\alpha } 
\notag \\
&&.E_{\beta }\left\{ -\frac{\beta }{1-\beta }(x-\xi )^{\beta +\alpha
}\right\} E_{\beta }\left\{ -\frac{\beta }{1-\beta }(t-\lambda )^{\beta
+\alpha }\right\} d\xi d\lambda .  \notag
\end{eqnarray}%
\noindent \noindent \textbf{Theorem 3 }: Let $f$ be a function such that $%
\frac{\partial ^{2}f}{\partial x\partial t}$ and $\frac{\partial ^{2}f}{%
\partial t\partial x}$ exist and continuous. Let $\alpha ,\beta \in (0,1)$
then the following equality holds.%
\begin{equation}
_{0}^{AC}D_{x,t}^{\alpha ,\beta }f(x,t)=_{0}^{AC}D_{t,x}^{\alpha ,\beta
}f(x,t)
\end{equation}

\noindent \noindent Proof : By definition we have that%
\begin{eqnarray}
&&_{0}^{AC}D_{x,t}^{\alpha ,\beta }f(x,t) \\
&=&\frac{A(\beta )}{1-\beta }\frac{1}{\Gamma \left\{ 1-\alpha \right\} }%
\int\limits_{0}^{x}\int\limits_{0}^{t}\frac{\partial ^{2}f(\xi ,\lambda )}{%
\partial \xi \partial \lambda }(x-\xi )^{-\alpha }(t-\lambda )^{-\alpha } 
\notag \\
&&.E_{\beta }\left\{ -\frac{\beta }{1-\beta }(x-\xi )^{\beta +\alpha
}\right\} E_{\beta }\left\{ -\frac{\beta }{1-\beta }(t-\lambda )^{\beta
+\alpha }\right\} d\xi d\lambda .  \notag
\end{eqnarray}%
Since $\frac{\partial ^{2}f}{\partial x\partial t}$ and $\frac{\partial ^{2}f%
}{\partial t\partial x}$ exist and are continuous then $\frac{\partial ^{2}f%
}{\partial x\partial t}=\frac{\partial ^{2}f}{\partial t\partial x}$
therefore we have%
\begin{eqnarray}
&&_{0}^{AC}D_{x,t}^{\alpha ,\beta }f(x,t) \\
&=&\frac{A(\beta )}{1-\beta }\frac{1}{\Gamma \left\{ 1-\alpha \right\} }%
\int\limits_{0}^{x}\int\limits_{0}^{t}\frac{\partial ^{2}f(\xi ,\lambda )}{%
\partial \xi \partial \lambda }(x-\xi )^{-\alpha }(t-\lambda )^{-\alpha } 
\notag \\
&&.E_{\beta }\left\{ -\frac{\beta }{1-\beta }(x-\xi )^{\beta +\alpha
}\right\} E_{\beta }\left\{ -\frac{\beta }{1-\beta }(t-\lambda )^{\beta
+\alpha }\right\} d\xi d\lambda .  \notag
\end{eqnarray}

\begin{eqnarray}
&=&\frac{A(\beta )}{1-\beta }\frac{1}{\Gamma \left\{ 1-\alpha \right\} }%
\int\limits_{0}^{t}\int\limits_{0}^{x}\frac{\partial ^{2}f(\xi ,\lambda )}{%
\partial \lambda \partial \xi }(t-\lambda )^{-\alpha }(x-\xi )^{-\alpha } 
\notag \\
&&.E_{\beta }\left\{ -\frac{\beta }{1-\beta }(t-\lambda )^{\beta +\alpha
}\right\} E_{\beta }\left\{ -\frac{\beta }{1-\beta }(x-\xi )^{\beta +\alpha
}\right\} d\lambda d\xi  \notag \\
&=&_{0}^{AC}D_{t,x}^{\alpha ,\beta }f(x,t)
\end{eqnarray}%
Thus 
\begin{equation}
_{0}^{AC}D_{x,t}^{\alpha ,\beta }f(x,t)=_{0}^{AC}D_{t,x}^{\alpha ,\beta
}f(x,t).
\end{equation}%
This completes the proof.

\section{Numerical Approximation}

There exist in nature many problems for which mathematical representations
are modeled with strong non-linearity. These kind of mathematical equations can mostly
handled using numerical approximations [8-10]. therefore in order to accommodate
researchers working in the field of numerical analysis. We present in this
section the numerical approximation of the $_{0}^{AC}D_{t}^{\alpha ,\beta
}f(t).$%
\begin{eqnarray}
&&_{0}^{AC}D_{t}^{\alpha ,\beta }f(t) \\
&=&\frac{A(\beta )}{1-\beta }\frac{1}{\Gamma \left\{ 1-\alpha \right\} }%
\int\limits_{0}^{t}\frac{df(\tau )}{d\tau }(t-\tau )^{-\alpha }E_{\beta
}\left\{ -\frac{\beta }{1-\beta }(t-\tau )^{\beta +\alpha }\right\} d\tau 
\notag \\
&=&\frac{A(\beta )}{1-\beta }\frac{1}{\Gamma \left\{ 1-\alpha \right\} }%
\int\limits_{0}^{t}\frac{f(t)-f(t+\Delta t)}{2\Delta t}(t-\tau )^{-\alpha
}E_{\beta }\left\{ -\frac{\beta }{1-\beta }(t-\tau )^{\beta +\alpha
}\right\} d\tau  \notag
\end{eqnarray}%
\begin{eqnarray*}
&&_{0}^{AC}D_{t}^{\alpha ,\beta }f(t) \\
&=&\frac{A(\beta )}{1-\beta }\frac{1}{\Gamma \left\{ 1-\alpha \right\} }%
\sum\limits_{i=0}^{n}\int\limits_{t_{i}}^{t_{i+1}}\frac{f^{i+1}-f^{i}}{%
2\left( \Delta t\right) }(t_{n}-\tau )^{-\alpha }E_{\beta }\left\{ -\frac{%
\beta }{1-\beta }(t_{n}-\tau )^{\beta +\alpha }\right\} d\tau \\
&=&\frac{A(\beta )}{1-\beta }\frac{1}{\Gamma \left\{ 1-\alpha \right\} }%
\sum\limits_{i=0}^{n}\frac{f^{i+1}-f^{i}}{2\left( \Delta t\right) }%
\int\limits_{t_{i}}^{t_{i+1}}(t_{n}-\tau )^{-\alpha }E_{\beta }\left\{ -%
\frac{\beta }{1-\beta }(t_{n}-\tau )^{\beta +\alpha }\right\} d\tau \\
&=&\frac{A(\beta )}{1-\beta }\frac{1}{\Gamma \left\{ 1-\alpha \right\} }%
\sum\limits_{i=0}^{n}\frac{f^{i+1}-f^{i}}{2\left( \Delta t\right) }%
\int\limits_{t_{n}-t_{i}}^{t_{n}-t_{i+1}}y^{-\alpha }E_{\beta }\left\{ -%
\frac{\beta }{1-\beta }y^{\beta +\alpha }\right\} dy \\
&=&\frac{A(\beta )}{1-\beta }\frac{1}{\Gamma \left\{ 1-\alpha \right\} }%
\sum\limits_{i=0}^{n}\frac{f^{i+1}-f^{i}}{2\left( \Delta t\right) }%
\int\limits_{t}^{x}y^{-\alpha }E_{\beta }\left\{ -\frac{\beta }{1-\beta }%
y^{\beta +\alpha }\right\} dy \\
&=&\frac{A(\beta )}{1-\beta }\frac{1}{\Gamma \left\{ 1-\alpha \right\} }%
\sum\limits_{i=0}^{n}\frac{f^{i+1}-f^{i}}{2\left( \Delta t\right) }\delta
_{n,i}.
\end{eqnarray*}%
Where%
\begin{eqnarray}
\delta _{n,i} &=&\left( t_{n}-t_{i+1}\right) ^{1-\alpha }E_{\beta ,2-\alpha
}\left\{ -\frac{\beta }{1-\beta }\left( t_{n}-t_{i+1}\right) ^{\beta +\alpha
}\right\} \\
&&-\left( t_{n}-t_{i}\right) ^{1-\alpha }E_{\beta ,2-\alpha }\left\{ -\frac{%
\beta }{1-\beta }\left( t_{n}-t_{i+1}\right) ^{\beta +\alpha }\right\} . 
\notag
\end{eqnarray}%
With the same time of idea we obtain the space approximation as 
\begin{equation}
_{0}^{AC}D_{x}^{\alpha ,\beta }f(x)=\frac{A(\beta )}{1-\beta }\frac{1}{%
\Gamma \left\{ 1-\alpha \right\} }\sum\limits_{i=0}^{n}\frac{f^{i+1}-f^{i}}{%
2\left( \Delta x\right) }\delta _{n,i}.
\end{equation}

\section{Conclusion}

In this paper we have broadened the scope of fractional calculus by
introducing fractional derivatives with two orders. The motivation of this
is from the fact that one can find in nature some systems with two different
and parallel layers with different properties, as for instance the movement
of water through aquifer with parallel layers with different properties.
This type of problem can not be portray with existing derivatives with
fractional order. Some properties connected to this new derivatives are
presented. The numerical approximation also presented. This new derivatives
will open new doors for PhD thesis, research and new fields.

\end{document}